\documentclass[12pt]{article}
\usepackage{epsfig}
\usepackage{bm}

\begin{document}

\title{Prime numbers: periodicity, chaos, noise}

\author{\small A. Bershadskii\\
\small ICAR, P.O. Box 31155, Jerusalem 91000, Israel }

\date{}
\maketitle

\begin{abstract}
Logarithmic gaps have been used in order to find a periodic component of the sequence 
of prime numbers, hidden by a random noise (stochastic or chaotic). The recovered period for 
the sequence of the first 10000 prime numbers is equal to $8\pm1$ (subject to the prime number theorem). 
For small and moderate values of the prime numbers (first 2000 prime numbers) this result has been 
directly checked using the twin prime killing method. 

\end{abstract}

\maketitle

~~~~~~~~~~~~~~~~~~~~~~~~~~~~~~~~~~~~~~~~~~~~~~~~~~~~~~~~~~~~~~~~~~~~~~~
\begin{flushright}
"There are some mysteries that      \\
the human mind will never penetrate.\\ 
To convince ourselves we have only  \\  
to case a glance at the tables of   \\ 
primes and we should perceive that  \\
there reigns neither order nor rule."\\

L. Euler
\end{flushright}

\section{Introduction}

The prime number distribution is apparently random. The apparent randomness 
can be stochastic or chaotic (deterministic). It is well known that the Riemann zeros obey the
chaotic GUE statistics (if the primes are interpreted as the classical periodic orbits of a 
chaotic system, see for a recent review \cite{bogom}), whereas the primes themselves 
are believed to be stochastically distributed (Poissonian-like etc. \cite{sou}). 
However, recent investigations suggest that primes themselves 
"...could be eigenvalues of a quantum system whose classical counterpart is chaotic at low energies 
but increasingly regular at higher energies." \cite{tt}. 
Therefore, the problem of chaotic (deterministic)
behavior of the moderate and small prime numbers is still open. Moreover, there are also 
an indication of periodic patterns in the prime numbers distribution \cite{tt}-\cite{ac}. These patterns, 
however, has been observed in {\it probability} distribution of the gaps between neighboring primes and 
not in the prime numbers sequence itself (see also Ref. \cite{dah}). The intrinsic randomness 
(stochastic or chaotic) of the primes distribution 
makes the problem of finding the periodic patterns in the prime numbers themselves 
(if they exist after all) a very difficult one.

\section{Logarithmic gaps}

It is believed that properties of the gaps between consecutive primes can provide a lot of information about 
the primes distribution in the natural sequence. The so-called prime number theorem states that 
the "average length" of the gap between a prime $p$ and the next prime number is proportional (asymptotically) 
to $\ln p $ (see, for instance, Ref. \cite{sou}). This implies, in particular, statistical non-stationarity of the prime numbers sequence. That is 
a serious obstacle for practical applications of the statistical methods to this sequence. For relatively 
large prime numbers one can try to overcome this obstacle by using relatively short intervals \cite{sou}. 
Another (additional) way to overcome this obstacle is to use logarithmic gaps (logarithms of the gaps). 
The non-stationarity in the sequence of the logarithmic gaps is considerably 'slower' than that in the original sequence of the gaps themselves. The logarithmic gaps have also another crucial advantage. If the noise 
corrupting the gaps sequence has a multiplicative nature, then for the logarithmic gaps this noise 
will be transformed into an additive one. It is well known that the additive noises can be much readily 
separated from the signal then the multiplicative ones (see below).

 Before starting the analysis let us recall two rather trivial properties of the gaps, which will be used below. 
The sequence can be restored from the gaps by taking cumulative sum. Pure periodicity in a sequence 
corresponds to a constant value of the gaps (the period). 

Let us take, as a first step, logarithms of the gaps for the sequence of prime numbers. Then, let us compute 
cumulative sum of these logarithms as a second step. 
Then, we will multiply each value in the cumulative sum by 10 and will replace each of the obtained values 
by a natural number which is nearest to it. As a result we will obtain a sequence of natural numbers: 
7,14,28,35,49,55,.... which will be called as ln-sequence.

Let us now define a binary function $v(n)$ of natural numbers $n=2,3,4,...$, which 
takes two values +1 or -1 and changes its sign passing any number from the ln-sequence. This function contains 
full information about the ln-sequence numbers distribution in the sequence of natural numbers.\\ 

\begin{figure} \vspace{-2cm}\centering
\epsfig{width=.7\textwidth,file=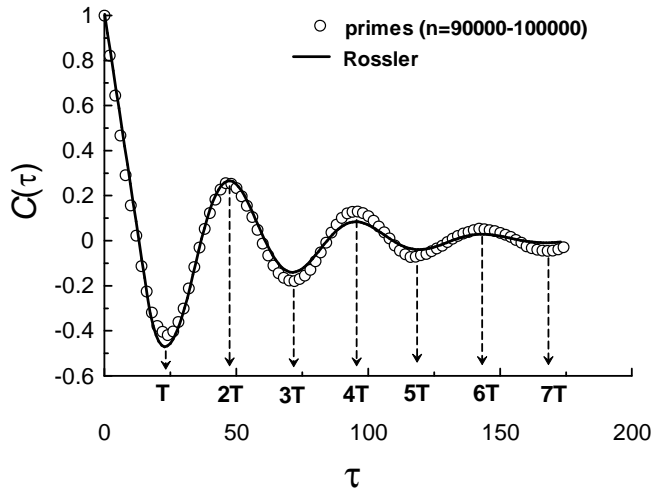} \vspace{-6.5cm}
\caption{Autocorrelation of the $v(n)$ functions for the ln-sequence (circles) in the interval $9\times 10^4 < n < 10^5$. The solid line correspond to autocorrelation function computed for the telegraph signal 
generated by the R\"{o}ssler attractor fluctuations overcoming the threshold $x = 7$. In order to make the autocorrelation functions comparable a rescaling has been made for the R\"{o}ssler autocorrelation function. }
\end{figure}
 
 Let us consider a relatively short interval: $9\times 10^4 < n < 10^5$. For this interval the autocorrelation 
function 
$$
C(n, \tau)= \langle v(n)v(n + \tau)\rangle - \langle v(n)\rangle \langle v(n + \tau)\rangle  \eqno{(2.1)}
$$
computed for the $v(n)$ function of the ln-sequence will be approximately independent on $n$. Figure 1 
shows the autocorrelation function (circles) computed for the $v(n)$ function of 
the ln-sequence for the 'short interval': $9\times 10^4 < n < 10^5$. Figure 2 shows an  autocorrelation function (circles) computed for the $v(n)$ function of the ln-sequence for interval of much smaller natural numbers: $1000 < n < 2\times 10^4$ (here we rely on the much slower 
non-stationarity in the logarithmic gaps sequence in comparison with the original gaps sequence). 
\begin{figure} \vspace{-2cm}\centering
\epsfig{width=.7\textwidth,file=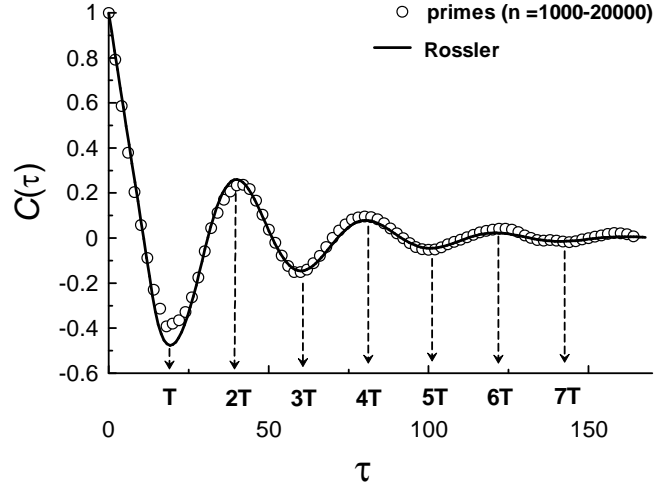} \vspace{-6.5cm}
\caption{As in Fig. 1 but for interval $1000 < n < 2\times 10^4$.
 }
\end{figure}

\section{Models}
\begin{figure} \vspace{-2cm}\centering
\epsfig{width=.7\textwidth,file=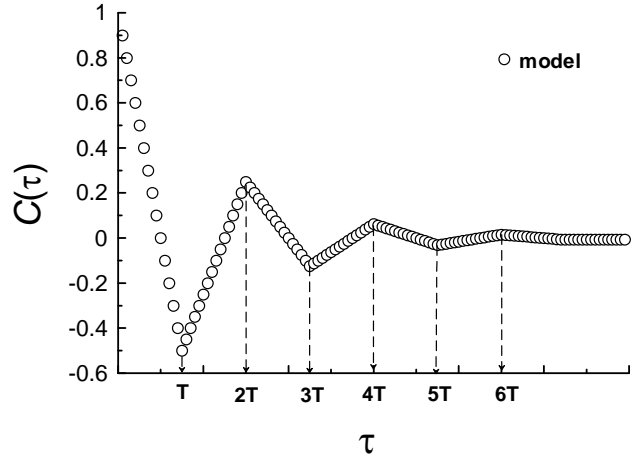} \vspace{-6.5cm}
\caption{Autocorrelation function for the simple model telegraph signal: Eq. (3.2) with $q=0.25$. 
 }
\end{figure}
In order to understand origin of the oscillating autocorrelation function shown in Figs. 1 and 2
let us consider a very simple telegraph signal, which allows analytic calculation of its 
autocorrelation function. The telegraph signal takes two values +1 or -1 and it changes its sign at 
discrete moments: 
$$
t_n=nT+\zeta \eqno{(3.1)}
$$
where $\zeta$ is an uniformly distributed over the interval $[0, T]$ random variable, 
$n=1,2,3...$ and $T$ is a fixed period. If $q$ 
is a probability of a sign change at a current moment ($0 \leq q < 1$), then the autocorrelation 
function of such telegraph signal is:
$$
C(\tau) = (n-\tau/T)(2q-1)^{n-1}+(\tau/T-(n-1))(2q-1)^n  \eqno{(3.2)}
$$
in the interval $(n-1)T \leq \tau < nT$. Figure 3 shows the autocorrelation function Eq. (3.2) calculated 
for $q=0.25$, as an example (cf. with Figs. 1 and 2). When parameter $q$ approaches value 0.5 the magnitude of 
the correlation function oscillations decreases . At value $q=0.5$ we have only linear decay in the interval 
$0 \leq \tau < T$ and for $\tau \geq T$ the correlation function $C(\tau) = 0$.\\
\begin{figure} \vspace{-2cm}\centering
\epsfig{width=.7\textwidth,file=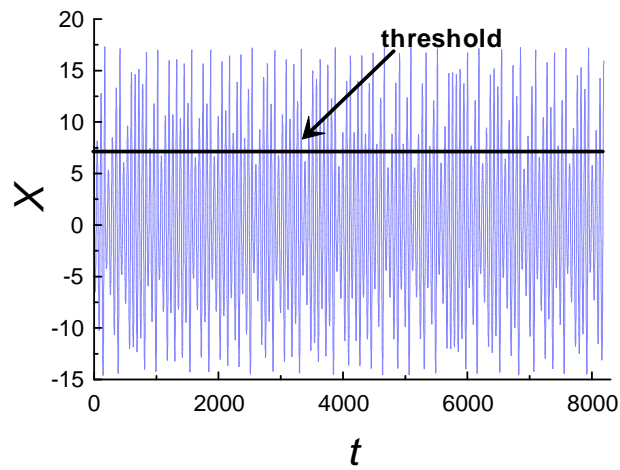} \vspace{-6.5cm}
\caption{X-component fluctuations 
of a chaotic solution of the R\"{o}ssler system Eq. (3.3) ($a=0.15,~ b=0.20,~ c=10.0$). }
\end{figure}

In the model Eq. (3.1) $\zeta$ was taken as a pure stochastic variable. This variable, however, can be also a chaotic (deterministic) one. Let us consider a chaotic solution of the R\"{o}ssler system \cite{ros} 
$$
\frac{dx}{dt} = -(y + z);~~  \frac{dy}{dt} = x + a y;~~  \frac{dz}{dt} = b + x z - c z\eqno{(3.3)}
$$      
where a, b and c are parameters (a reason for this system relevance will be given in Section 4). 
Figure 4 shows the x-component fluctuations of a chaotic solution 
of the R\"{o}ssler system.  

Let us consider a telegraph signal generated by the variable $x$ crossing certain threshold from below (Fig. 4). 
This telegraph signal takes values: +1 or -1, and changes its sign at the threshold-crossing points. 
Then, the fundamental period of the R\"{o}ssler chaotic attractor provides the period $T$ in Eq. (3.1), 
whereas the chaotic fluctuations of the variable $x$ of the attractor provide chaotic variable 
$\zeta$ in the Eq. (3.1) for the considered telegraph signal. Autocorrelation function computed for this telegraph 
signal has been shown in Figs. 1 and 2 as the solid line. In order to make the autocorrelation functions comparable a rescaling has been made for the R\"{o}ssler telegraph signal's autocorrelation function.\\ 

Now, let us restore the fundamental period $T_0$ for the sequence of the first 10000 prime numbers:
$$
T_0= \exp (T/10) \simeq 8\pm 1 \eqno{(3.4)} 
$$
where the period $T$ for the ln-sequence was taken from Figs. 1 and 2. 

\section{Twin prime killing method}
\begin{figure} \vspace{-2cm}\centering
\epsfig{width=.7\textwidth,file=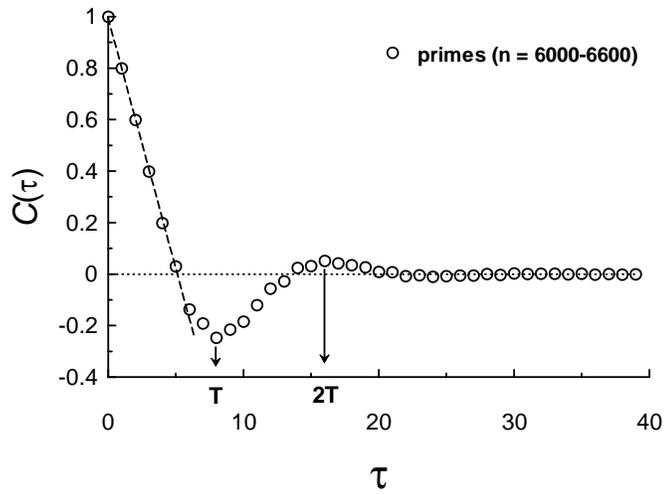} \vspace{-6.5cm}
\caption{Autocorrelation of the $v(n)$ function for the K2-subsequence (circles) in the relatively short 
interval $6000 < n < 6600$.  }
\end{figure}
One can expect that the smallest gaps between prime numbers are a source of 
the non-additive (high frequency) noise. Especially  for small and moderate prime numbers 
where the minimal gap $\Delta = 2$ occurs rather frequently. The primes separated by such gap 
are called twin primes. 
Let us kill one (say larger one) of the twin primes from each pair of the twin primes (a smart high frequency filter). 
Due to remaining part of the twins the subsequence of the primes preserves roughly the structure of the original full sequence of the primes (let us call this subsequence as K2-subsequence). 

\begin{figure} \vspace{-2cm}\centering
\epsfig{width=.7\textwidth,file=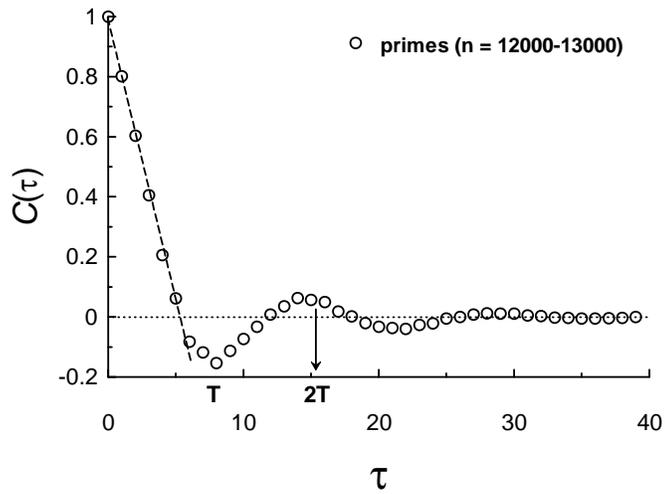} \vspace{-6.5cm}
\caption{As in Fig. 5 but for interval $12000 < n < 13000$.  }
\end{figure}
\begin{figure} \vspace{-1.5cm}\centering
\epsfig{width=.7\textwidth,file=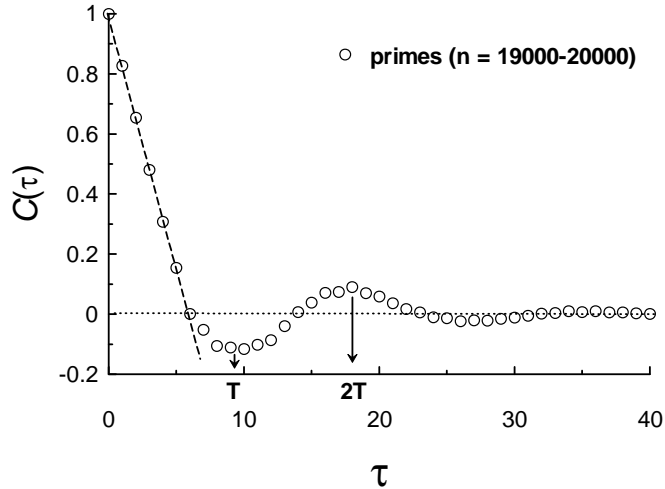} \vspace{-6.5cm}
\caption{As in Fig. 5 but for interval $19000 < n < 20000$.  }
\end{figure}
On the other hand, for small and moderate prime numbers the K2-subsequence can be expected to be less 
corrupted by the non-additive high frequency noise than the original sequence of the primes. Indeed, figure 5 shows 
autocorrelation of the $v(n)$ function for the K2-subsequence (circles) in the relatively short 
interval $6000 < n < 6600$. The dashed straight line indicates the characteristic linear decay (cf. with Figs. 1-3)  
and the arrows indicate the period $T = 8 \pm 1$. Figures 6 and 7 shows analogous results for the intervals $12000 < n < 13000$ and $19000 < n <20000$ respectively. Comparing this result with Eq. (3.4) 
we can see good agreement between the values of the hidden period of the prime number sequence obtained by the 
two different methods: logarithmic gaps and twin prime killing. 

  It should be noted that the most probable value of the gap is $\Delta =6$ for the both original and 
K2-subsequence (one of the two integers nearest to any prime: $p -1$ or $p+1$, always has 6 as its divisor). 
Although the hidden period does not coincide with the most probable value of the gap (and apparently it is a 
subject of the prime number theorem), the simple 6-divisor could be a starting point in an attempt to understand 
origin of the hidden periodicity (in particular, the end-point of the linear decay in the Figs. 5-7 is always equal to 6, cf. with Fig. 3).


\newpage
\begin{thebibliography}{99}
\bibitem{bogom} E. Bogomolny, Riemann zeta function and quantum chaos, 
Prog. Theor. Phys. Supplement, vol. 166,  pp. 19-44, 2007.
\bibitem{sou} K. Soundararajan, The distribution of prime numbers, NATO Science Series II: 
Mathematics, Physics and Chemistry, , vol. 237, pp. 59-83, 2007.
\bibitem{tt} T. Timberlake and J.Tucker, Is there quantum chaos in the prime numbers?, 
Bulletin of the American Physical Society,  vol. 52, p. 35, 2007, (arXiv:0708.2567).
\bibitem{po} C.E. Porter, Statistical Theories of Spectra: Fluctuations, NY: Academic Press, 1965.
\bibitem{wo} M. Wolf, Applications of statistical mechanics in
number theory, Physica A, vol. 274, pp. 149-157, 1999.
\bibitem{kis} P. Kumar, P.Ch. Ivanov, and H.E. Stanley, Information Entropy and Correlations in Prime Numbers, 
arXiv:cond-mat/0303110.
\bibitem{ac} S. Ares, and M. Castro, Hidden structure in the randomness
of the prime number sequence?, Physica A, vol. 360, pp. 285-296, 2006.
\bibitem{dah} S. R. Dahmen, S. D. Prado and T. Stuermer-Daitx, Similarity in the statistics of prime numbers, 
Physica A, vol. 296, pp. 523-528, 2001. 
\bibitem{ros} O.E. R\"{o}ssler, An equation for continuous chaos, Phys. Lett. A,  vol. 57, pp. 397-398, 1976.

\end{thebibliography}
\end{document}